\documentclass{commat}

\newcommand{\ds}{}
\newcommand{\row}[2]{#1_1,\ldots,#1_{#2}}
\newcommand{\modd}[1]{\ \text{\rm (mod}\, #1\text{\rm )}}
\newcommand{\natop}[2]{\genfrac{}{}{0pt}{2}{#1}{#2}}
\newcommand{\Z}{\mathbb{Z}}
\newcommand{\ca}{\mathcal {A}}

\title{Congruences for the cycle indicator of the symmetric group}
\author{Abdelaziz Bellagh and Assia Oulebsir}

\affiliation{
    \address{Abdelaziz\, Bellagh (Corresponding author) -- Facult\'e de
Math\'ematiques, U.S.T.H.B, Laboratoire LATN, B.P 32, El ALIA, Bab
 Ezzouar, 16111, Alger, {\sc Alg\'erie}.}
  \email{abellagh@yahoo.fr}
   \address{Assia\, Oulebsir -- Facult\'e de
Math\'ematiques, U.S.T.H.B, Laboratoire LATN, B.P 32, El ALIA, Bab
Ezzouar, 16111, Alger, {\sc Alg\'erie}.}
  \email{oulebsir.assia@hotmail.com}
    }

\abstract{Let $n$ be a positive integer and let $C_n$ be the cycle indicator
of the symmetric group $S_n$. Carlitz proved that if $p$ is a
prime, and if $r$ is a non negative integer, then we have the
congruence $$C_{r+np}\equiv (X_1^p-X_p)^nC_r
 \modd{p\Z_p[X_1,\cdots,X_{r+np}]},$$
where $\Z_p$ is the ring of $p$-adic integers. 
We prove that for $p\neq 2$, the preceding congruence holds modulo 
$np\Z_p[X_1,\cdots,X_{r+np}]$. This allows us to
prove a Junod's conjecture for Meixner polynomials.}



\VOLUME{31}
\YEAR{2023}
\NUMBER{1}
\firstpage{251}
\DOI{https://doi.org/10.46298/cm.10391}

\begin{paper}

\section{Introduction and results}
Let $n$ be a positive integer. The cycle indicator of the symmetric
group $S_n$, is the polynomial $C_n$ defined by 
\begin{equation*}
C_n=\sum
c_n(\row{m}{n})\prod_{i=1}^nX_i^{m_i},
\end{equation*}
 where the sum is over all non
negative integers $m_i$ such that  $\sum_{i=1}^nim_i=n$, and 
\begin{equation*}
c_n(\row{m}{n})=\frac{n!}{\prod_{i=1}^n
  i^{m_i}(m_i!)}
\end{equation*}
If we set $C_0=1$, then the exponential generating function for $C_n$
is 
\begin{equation}
\ds \exp\Big(\sum_{i=1}^{\infty}X_i\frac{t^i}{i}\Big)=
\sum_{n=0}^{\infty}C_n(\row{X}{n})\frac{t^n}{n!};\label{foncgencn}
\end{equation}
furthermore, the coefficients $c_n(\row{m}{n})$ of the polynomial
$C_n$ are integers; see Riordan \cite[p.\ 67--68]{Key5}.

Carlitz \cite[p\ 1222]{Key1} proved, that for any prime $p$, 
and any positive integer $n$, and $m_1,\ldots,m_{np}\in \{0,\ldots,np\}$, such
that $np=\sum_{i=1}^{np}im_i$, and for any non negative integer $r$, we have
\begin{equation}
c_{np}(\row{m}{np})\equiv \begin{cases}
\ds (-1)^{m_p}\dbinom{n}{m_p}\modd{p\Z_p},& \text{ if}\;
\sum_{i\notin\{1,p\}}m_i=0,\\
\hfill\\ 
0 \modd{p\Z_p},               & \text{ if}\;
\sum_{i\notin\{1,p\}}m_i\neq0;\label{eq00}\end{cases}
\end{equation}
\begin{equation}
\text{and}\;\; C_{r+np}\equiv (X_1^p-X_p)^nC_r
\modd{p\Z_p[X_1,\cdots,X_{r+np}]},
\label{eq01}
\end{equation}
where $\Z_p$ is the ring of $p$-adic integers.
We give the following generalization:
\begin{proposition}
Let $r, n, p, m_1,\ldots,m_{np}$ be as above, and let 
\begin{equation*}
n^*=\begin{cases} n/2,  \, \mbox{\rm if}\; n\in
  2\Z\\
n,  \, \mbox{\rm otherwise}.\\
\end{cases}
\end{equation*}
1) If $\;\sum_{i\notin\{1,p\}}m_i=0$, then 
\begin{equation}
c_{np}(\row{m}{np})\equiv 
\ds (-1)^{pm_p}\dbinom{n}{m_p}\modd{n^*p\Z_p} \label{congr1p}.
\end{equation}
2) If $\;\sum_{i\notin\{1,p\}}m_i\neq 0$, then
\begin{equation}  
 c_{np}(\row{m}{np})\equiv \;0 \modd{n^*p\Z_p} \label{congr2p}.            
\end{equation}
3) The cycle indicator polynomials satisfy the congruence 
\begin{equation}
C_{r+np}\equiv (X_1^p-X_p)^nC_r
\modd{n^*p\Z_p[X_1,\cdots,X_{r+np}]}. \label{congr3p}
\end{equation}
\end{proposition} 
\begin{corollary}\label{complementary}
Let $p$ be a prime, and let $n, r$ be positive integers such
that $1\leq r\leq p-1$ and $r+np=\sum_{i=1}^{r+np}im_i$, where the
$m_i$ are non negative integers. We define $n^*$
as in the proposition. 
\begin{enumerate}
\item If $m_1\geq p(n-m_p)\geq 0$, then
we have the following congruence modulo $n^*p\Z_p$,
\begin{equation}
 \hskip-0.6cm c_{r+np}(\row{m}{r+np})\equiv (-1)^{pm_p}\dbinom{n}{m_p}
c_{r}(m_1+pm_p-np,m_2,\ldots,m_r).\label{congr1co}
\end{equation}
\item Otherwise, we have $ c_{r+np}(\row{m}{r+np})\equiv 0
   \modd{n^*p\Z_p}$.
\end{enumerate}
\end{corollary}
\begin{remark}
Under the hypotheses of Corollary \ref{complementary}, if
$r+pn=m_1+pm_p$ and if $n\geq m_p$, 
then we deduce from the congruences \eqref{congr1p} and \eqref{congr1co}, that 
\begin{equation*}
c_{r+np}(m_1,0,\ldots 0,m_p,0\ldots,0)\equiv c_{np}(m_1-r,0,\ldots 0,m_p,0\ldots,0) \modd{n^*p\Z_p}
\end{equation*}
\end{remark}
 We will need the following lemma in the next corollary.
\begin{lemma}[Junod \cite{Key3}]
 Let $m,n$ be two  positive integers, and let $\alpha, \beta$ be two
 elements of a commutative ring
 $\ca$ containing $\Z_p$. 

If $m\in p\Z$ and $\alpha\equiv \beta \modd{m \ca}$, 
then $\alpha^n\equiv \beta^n \modd{mn \ca}$. 
\end{lemma}

Meixner polynomials are defined by their exponential generating
function 
\begin{equation*}
\ds \frac{1}{\sqrt{1+t^2}}\exp (X \arctan
t)=\sum_{n=0}^{\infty}Q_n(X)\frac{t^n}{n!}.
\end{equation*}
Let $Q_n^{\ast}$ be the polynomials defined by the exponential generating
function 
\begin{equation*}
\ds \exp (X \arctan
t)=\sum_{n=0}^{\infty}Q_n^{\ast}(X)\frac{t^n}{n!}.
\end{equation*}
 Junod \cite[p\ 73]{Key2}, proved that, if $p\neq 2$, then
\begin{gather}
Q_{np}^{\ast}(X)\equiv Q_{np}(X) \modd{np\Z_p[X]},\;\text{and}\label{congrj1}\\
Q_{p}(X)\equiv (X^p-(-1)^{(p-1)/2}X) \modd{p\Z_p[X]}.\label{congrj2}
\end{gather}
Since $X \arctan t=\sum_{i=1}^{\infty}x_i\frac{t^i}{i}$, where $x_i = 0$, when $i$ is even, and $x_i = (-1)^{(i-1)/2}X$, when $i$ is odd, it follows from equality
\eqref{foncgencn}, that $Q_n^{\ast}(X)=C_n(\row{x}{n})$, for 
any positive integer $n$.
Using \eqref{congr3p}, \eqref{congrj1},
\eqref{congrj2} and the lemma, we deduce that:
\begin{corollary}
If $p\neq 2$, then 
\begin{equation*}
Q_{np}(X)\equiv Q_{p}^n(X)\equiv (X^p-(-1)^{(p-1)/2}X)^n
\modd{np\Z_p[X]}.
\end{equation*}
\end{corollary}

This gives us a positive answer to a conjecture of Junod \cite[p\
 74]{Key2}.
\section{Proof of the proposition}
For $n\notin p\Z$, the congruences \eqref{congr1p} and \eqref{congr2p}
follow from \eqref{eq00}, or from Macdonald \cite[p\ 30]{Key4}.
Hence, we prove these congruences for $n\in p\Z$.

1) Let $m$ be a non negative integer, then we have
\begin{equation}
\dbinom{np}{pm}\equiv \dbinom{n}{m} \modd{np\Z_p},\;\,\text{and}\;\; pm\dbinom{n}{m}\equiv 0 \modd{np\Z_p}\label{congr1l1}\\
\end{equation}
\begin{gather}
(mp)!=(-1)^{pm+1}\Gamma(pm+1)(m!)p^{m},\;\,\text{and}\,\label{eq2l1}\\
\Gamma(pm+1)+1 \in (pm/2)\Z_p\,\label{congr3l1},
\end{gather}
where $\Gamma$ denotes the Morita $p$-adic Gamma function;  
see Robert \cite[p 369]{Key4}. 

If $\sum_{i\notin\{1,p\}}m_i=0$, then we have
\begin{equation*}
\ds c_{np}(\row{m}{np})=\frac{(np)!}{m_1!m_p!p^{m_p}}=
\dbinom{np}{pm_p}\frac{(pm_p)!}{{m_p}!p^{m_p}}.
\end{equation*}
Hence the congruence \eqref{congr1p} follows from \eqref{congr1l1}, 
\eqref{eq2l1} and \eqref{congr3l1}.
Furthermore, in this case, $c_{np}(\row{m}{np})$ and $\dbinom{n}{m_p}$ have the same $p$-adic valuation, since
\begin{equation}
\ds c_{np}(\row{m}{np})=
(-1)^{pm_p}\dbinom{n}{m_p}\frac{\Gamma(np+1)}{\Gamma(m_1+1)}\label{formula2l1}.
\end{equation}

2) Now we prove by induction on the
$p$-adic valuation $\nu$ of $n$, that if $\ds
np=\sum_{i=1}^{np} im_i$ and if there exists
$i\in \{2,\ldots,np\}-\{p\}$ such that $m_i\neq 0$, then 
the congruence \eqref{congr2p} is satisfied. For $\nu =0$, the congruence
\eqref{congr2p} is true.
We assume the congruence \eqref{congr2p} holds for $\nu-1\geq 0$.

 First, we consider the case where $p$ does not divide $m_i$. Then we have
\begin{gather*}
\ds np-i=i(m_i-1)+\sum_{\natop{j=1}{j\neq i}}^{np-i}jm_j,\;\text{and}\;\\
\ds
c_{np}(m_1,\ldots,m_{np})=\Big(\frac{np(np-1)\cdots(np-i+1)}{im_i}\Big)
c',\,
\text{where}\\
 c'= \begin{cases}
c_{np-i}(m_1,\ldots,m_{i-1},m_i-1,m_{i+1},\ldots,
m_{np-i}), \text{ if}\; np\geq 2i ,\\
\hfill\\ 
c_{np-i}(m_1,\ldots, m_{np-i}), \text{ if}\;  np< 2i \;\text{and}\; m_i=1 .\label{eq000}\end{cases}
\end{gather*}
Then, we note that $(np-1)\cdots (np-i+1)\in (i/2)\Z_p $.

Indeed, if $\mu$ is
 the $p$-adic valuation of $i$, and if $\mu\geq 2$ with
 $i\neq 4$, we have
 $$(np-1)\cdots (np-i+1)\in (np-p^{\mu-1}-p)\ds
 \prod_{j=1}^{\mu-1}(np-p^j)\Z_p .$$
 Hence the congruence 
 \eqref{congr2p} follows (in this case, we do not need the induction hypothesis).

On the other hand, if for all $j\notin \{1,p\}$, we have that $p$ 
divides $m_j$, then
$p$ divides $m_1$. Hence, if we set $m_j=pm'_j$ for any $j\neq p$, 
we obtain that
\begin{equation*} 
\ds n=(m'_1+m_p)+\sum_{\natop{j=2}{j\neq p}}^{n}jm'_j.
\end{equation*}
By \eqref{eq2l1}, we get
\begin{gather*}
 c_{np}(m_1,\ldots,m_{np})=uz\dbinom{m_p+m'_1}{m_p} c',\,\text{where}\,\\
\ds u=(-1)^{np+1}\Gamma(pn+1)\prod_{\natop{j=1}{j\neq
    p}}^{n}\Gamma(m'_jp+1)^{-1}(-1)^{pm'_j+1},\\
z=\prod_{\natop{j=2}{j\neq p}}^{n}p^{(j-1)m'_j}j^{-m'_j(p-1)},\;
\text{and}\,\\
 c'= c_{n}(m'_1+m_p,m'_2,\ldots,m'_{p-1},0,m'_{p+1},\ldots,m'_{n}).
\end{gather*}
Since the $p$-adic Gamma function takes its values in the set of
inversible elements in the ring $\Z_p$, we deduce that the $p$-adic
valuation of $u$ is zero. Then we remark that if $\mu$ denotes the $p$-adic
valuation of a positive integer $i$ such that $i\notin\{1, p\}$, we have $i-1>(p-1)\mu$. Hence $z\in
p\Z_p$. By the induction hypothesis (on the
$p$-adic valuation $\nu$ of $n$), we have $c'\in (n/p)^*p\Z_p$. Since $(n/p)^*p\in
n^*\Z_p$, we conclude that 
the congruence \eqref{congr2p} holds for $\nu$.

3) Using the congruences \eqref{congr1p} and \eqref{congr2p}, we get 
\begin{align}
C_{np} &\equiv 
 \sum_{m_1+pm_p=np}(-1)^{pm_p}\dbinom{n}{m_p}X_1^{m_1}X_p^{m_p} 
\modd{n^*p\Z_p[X_1,\cdots,X_{np}]},\notag \\
           &\equiv \;(X_1^p+(-1)^pX_p)^n \modd{n^*p\Z_p[X_1,\cdots,X_{np}]}. \label{congr4}
\end{align}
In particular, if $\nu$ is the $p$-adic valuation of $n$, we have
\begin{equation}
C_{p^{\nu+1}}\equiv (X_1^p+(-1)^pX_p)^{p^{\nu}} \modd{n^*p\Z_p[X_1,\cdots,X_{np}]}. \label{congr5}
\end{equation}
Deriving equality \eqref{foncgencn} with respect to $t$, we obtain
that for any positive integer $m$, we have
\begin{equation}
C_m=\ds \sum_{j=0}^{m-1} \frac{(m-1)!}{j!}X_{m-j}C_j=
\begin{array}{|ccccc|}
X_1    &     -1  & 0      &   & 0 \\
X_2    &    X_1  & -2     &   & 0 \\
\vdots &  \vdots & \vdots &   & \vdots\\
X_{m-1}& X_{m-2} & X_{m-3}&   & -(m-1)\\
X_m    & X_{m-1} & X_{m-2}&   &  X_1 \\
\end{array} \label{eq1}
\end{equation}
(by expanding the determinant by the last row).

Taking $k=n/p^{\nu} $, $m=r+np$, and 
reducing the identity \eqref{eq1} modulo $p^{\nu+1}$, we get
\begin{equation*}
 C_{r+kp^{\nu+1}}\equiv C_r(C_{p^{\nu+1}})^{k} \modd{p^{\nu+1}\Z_p[X_1, \dotsc, X_m]}.
\end{equation*}
 Using \eqref{congr5}, we deduce
the congruence \eqref{congr3p}.\qed
\begin{remark}
If $p\neq 2$, $(p-1)!\not\equiv -1\modd{p^2\Z}$ (i.e., if $p$ is not a Wilson's
prime), and if $n\in p\Z$, the congruence \eqref{congr1p}, does not hold modulo
$np^2\Z_p$. 
Indeed, by \eqref{formula2l1}, we have for $m_p=1$,
\begin{equation*}
c_{np}(\row{m}{np})-\ds (-1)^{pm_p}\dbinom{n}{m_p}=n (1-q).
\end{equation*}
where
$q=\dfrac{\Gamma(np+1)}{\Gamma(m_1+1)}=-\prod_{j=1}^{p-1}(np-j)\equiv 
-\prod_{j=1}^{p-1}j \modd{p^2\Z}$.
\end{remark}


\EditInfo{June 19, 2021}{October 08, 2021}{Pasha Zusmanovich}

\end{paper}
\begin{references}

\refer{Paper}{Key1}
\Rauthor{Carlitz. L}
\Rtitle{Some congruences for the Bell polynomials}
\Rjournal{Pacific Math. J}
\Rvolume{Vol. 11}
\Ryear{1961}
\Rnumber{No. 4}
\Rpages{1215-1222}

\refer{Other}{Key2}
\Rauthor{Junod.A}
\Rtitle{Congruences par l'analyse $p$-adique
et le calcul symbolique. Th\`ese de Doctorat Universit\'e de Neuch\^atel.}
\Ryear{2003}

\refer{Paper}{Key3}
\Rauthor{Junod.A}
\Rtitle{ Congruences pour les polyn\^omes et nombres de Bell}
\Rjournal{Bulletin de la Soci\'et\'e Math\'ematique de Belgique}
\Rvolume{ }
\Ryear{2002}
\Rnumber{No. 9}
\Rpages{503-509}

\refer{Book}{Key4}
\Rauthor{Macdonald.I.G}
\Rtitle{Symmetric function and Hall polynomials}
\Rpublisher{Oxford science publications mathematical monographs}
\Ryear{1995}

\refer{Book}{Key5}
\Rauthor{Riordan.J}
\Rtitle{An introduction to combinatorial analysis}
\Rpublisher{John Wiley and Sons, Inc, New York, London, Sydney}
\Ryear{1967}

\refer{Book}{Key6}
\Rauthor{Robert.A}
\Rtitle{A course in p-adic analysis}
\Rpublisher{Springer-Verlag, New York-Heidelberg}
\Ryear{1999}

\end{references}
